\documentclass[12pt]{article}
\usepackage[top=2.25cm, bottom=2.25cm, left=2.25cm, right=2.25cm]{geometry}
\usepackage[T1]{fontenc}
\usepackage{lmodern}
\usepackage{amsmath,amssymb,amsthm,mathtools}
\usepackage{microtype,booktabs,array}
\usepackage[hidelinks]{hyperref}
\allowdisplaybreaks[2]
\newtheorem{theorem}{Theorem}[section]
\newtheorem{lemma}[theorem]{Lemma}

\newtheorem{corollary}[theorem]{Corollary}
\newtheorem{claim}[theorem]{Claim}
\newtheorem{conjecture}[theorem]{Conjecture}
\theoremstyle{remark}
\newtheorem*{remark}{\bf Remark}

\title{Judicious $r$-partitions of $r$-uniform hypergraphs}
\author{Siwei Lin\footnote{Email: linsw0710@163.com},~~ Qinghou Zeng\footnote{Research partially supported by National Key R\&D Program of China (Grant No. 2023YFA1010202) and National Natural Science Foundation of China (Grant No. 12371342). Email: zengqh@fzu.edu.cn (Corresponding author)}\\
{\small Center for Discrete Mathematics, Fuzhou University, Fujian, 350003, China}}
\date{}
\begin{document}
\maketitle

\begin{abstract}
Bollob\'as and Thomason conjectured that every $r$-uniform hypergraph with $m$ edges admits an $r$-partition such that each class meets at least $rm/(2r-1)$ edges.
In this paper, we establish a lower bound of $(11-\sqrt{13})m/18$ for every $r\ge 2$. This improves Haslegrave's bound $rm/(3r-4)$ for every $r\ge 8$.
\end{abstract}

\section{Introduction}

Many classical graph and hypergraph partitioning problems seek to optimize a single quantity. A fundamental example is the Max-Cut problem, which asks for a bipartition $V(G)=V_1\cup V_2$ maximizing the number $e(V_1,V_2)$ of crossing edges. Edwards \cite{Edwards73,Edwards75} proved the essentially sharp bound $e(V_1,V_2)\ge m/2+(\sqrt{8m+1}-1)/8$ for every graph with $m$ edges. Bollob\'as and Scott \cite{BS2002} extended this result to $k$-partitions, showing that every graph with $m$ edges admits a partition $V(G)=V_1\cup\cdots\cup V_k$ such that $e(V_1,\ldots,V_k) \ge \frac{k-1}{k}(m+(\sqrt{8m+1}-1)/4)-O(k)$, where $e(V_1,\ldots,V_k)$ denotes the number of edges joining distinct classes.

In practice, however, one often seeks a partition of a graph or hypergraph that optimizes several quantities simultaneously. Bollob\'as and Scott called such problems judicious partitioning problems. 
A basic example is the bottleneck bipartition problem, posed by Entringer \cite{SS1994}, which asks for a bipartition $V(G)=V_1\cup V_2$ minimizing $\max\{e(V_1),e(V_2)\}$, where $e(V_i)$ denotes the number of edges induced by $V_i$. 
Bollob\'as and Scott \cite{BS1999} proved that every graph with $m$ edges admits a bipartition such that, for $i=1,2$, $e(V_i)\le \frac m4+\frac{\sqrt{8m+1}-1}{16}$.
This bound is sharp for complete graphs of odd order. They also considered the corresponding problem for hypergraphs \cite{BS1997}, proving that every $3$-uniform hypergraph with $m$ edges admits a bipartition in which each class contains at most $m/8+5m^{6/7}$ edges. 
More generally, in the same paper, they conjectured that every $r$-uniform hypergraph with $m$ edges admits a bipartition such that each class contains at most $m/2^r+o(m)$ edges. This conjecture remains open in general.

Another important direction in judicious partitioning is to find a partition in which every class meets many edges. For a hypergraph $H$ and $X\subseteq V(H)$, let $d_H(X)$ denote the number of edges meeting $X$, or their total weight when $H$ is weighted. A particularly interesting case occurs when we partition the vertices of an $r$-uniform hypergraph into $r$ classes. Bollob\'as and Thomason (see \cite{BRT1993,BS2000}) proposed the following conjecture.

\begin{conjecture}[Bollob\'as and Thomason (see \cite{BRT1993,BS2000})]\label{BT-Conj}
Every $r$-uniform hypergraph with $m$ edges admits a partition $V(H)=V_1\cup\cdots\cup V_r$ such that, for $1\le i\le r$,
\[
d_H(V_i)\ge \frac{r}{2r-1}m.
\]
\end{conjecture}
For \(r=2\), the conjecture follows directly from the local optimality of a maximum cut. For \(r=3\), Bollob\'as and Scott \cite{BS2000} proved that every \(3\)-uniform hypergraph with \(m\) edges admits a tripartition in which each part meets at least \(\frac{5m-1}{9}\) edges. Haslegrave \cite{Haslegrave2012} subsequently improved this bound to \(0.6m\), thereby confirming the conjecture for \(r=3\). Further asymptotic improvements were obtained by Ma and Yu \cite{MY2012}, who established a bound of \(0.65m+o(m)\), and by Spink and Tiba \cite{ST2020}, who obtained the optimal bound \(\frac{19}{27}m+O(m^{6/7})\).

For general \(r\ge 3\), Bollobás, Reed and Thomason \cite{BRT1993} proved that every \(r\)-uniform hypergraph with \(m\) edges admits a partition in which each part meets at least \(\frac{1-1/e}{3}m\) edges. Haslegrave \cite{Haslegrave2014} subsequently improved this bound to \(\frac{rm}{3r-4}\). More recently, Lin and Zeng \cite{LZ2026} further improved the bound to \(\frac{37}{64}m-o(m)\)  for \(r=4\), confirming Conjecture~\ref{BT-Conj} for $4$-uniform hypergraphs with sufficiently many edges.

For further results on judicious partitioning, see \cite{BS1993,BS1997,BS2010,HWY2016,HX2027,JMYY2019,LSSZ2019,MY2016,ST2020,XYY2010,XY2011,XY2014}.

Our main result provides a new bound, improving Haslegrave's estimate for every $r\ge8$. Moreover, it holds for nonnegatively weighted $r$-uniform hypergraphs. 

\begin{theorem}\label{thm:main}
For every integer $r\ge2$, every finite nonnegatively weighted $r$-uniform hypergraph $H$ of total weight $1$ admits a partition $V(H)=V_1\cup\cdots\cup V_r$ such that, for $1\le i\le r$,
\[
d_H(V_i)\ge \frac{11-\sqrt{13}}{18}.
\]
\end{theorem}

\begin{corollary}
For every integer $r\ge2$, every finite $r$-uniform hypergraph $H$ with $m$ edges admits a partition $V(H)=V_1\cup\cdots\cup V_r$ such that, for $1\le i\le r$,
\[
d_H(V_i)\ge
\frac{11-\sqrt{13}}{18}m.
\]
\end{corollary}
\begin{remark}
    Note that $\frac{11-\sqrt{13}}{18}\approx 0.41>\frac{r}{3r-4}$ when $r\ge8$.
    Thus, our bound improves Haslegrave's result for every $r\ge8$.
    We mention that our approach can be refined for small $r$.
    We omit the details, as the resulting improvements are modest.
\end{remark}

\section{Hypergraphs with few vertices of large degree}\label{sec:fractional}

A set function $f:2^U\to\mathbb R_{\ge0}$ is \emph{normalized} if $f(\varnothing)=0$; its multilinear extension is
\[
 F(x)=\mathbb E [f(R_x)],
\]
where the elements of $U$ are included in $R_x$ independently, with inclusion probabilities $x_u$.
A set function $f:2^U\to\mathbb R_{\ge0}$ is \emph{monotone} if $f(X)\le f(Y)$ for all $X\subseteq Y\subseteq U$, and it is \emph{submodular} if $f(X\cup\{u\})-f(X)\ge f(Y\cup\{u\})-f(Y)$ for $u\in U$ and $X\subseteq Y\subseteq U$.

\begin{lemma}\label{lem:round}
Let $f_1,\dots,f_k$ be normalized, monotone, submodular functions on a
finite set $U$, and let $F_1,\ldots,F_k$ be their multilinear extensions. If
$x_{iu}\ge0$ and $\sum_{i=1}^k x_{iu}=1$ for every $u\in U$, then $U$
has a partition $P_1\cup\cdots\cup P_k$ such that
\[
 f_i(P_i)\ge F_i(x_i)-D_i,
\]
where $D_i=\max_{u\in U}f_i(\{u\})$. For $U=\varnothing$, we take $D_i=0$.
\end{lemma}
\begin{proof}
Define the \textit{support graph} of $(x_{iu})$ to be the bipartite graph with vertex set $U\cup[k]$ and edge set $\{\{i,u\}:i\in[k],\,u\in U,\,x_{iu}>0\}$.
For a monotone submodular function, $\partial F_i/\partial x_{iu}\ge0$, and
$\partial^2F_i/(\partial x_{iu}\partial x_{iv})\le0$ when $u\ne v$.

Suppose the support graph contains a cycle $i_1u_1i_2u_2\cdots i_hu_hi_1$, with cyclic subscripts. 
If the derivative associated with an edge of the cycle, say $i_1u_1$, is zero, then transfer the entire mass from $i_1u_1$ to $i_2u_1$; that is, set $x'_{i_1u_1}=0$ and $x'_{i_2u_1}=x_{i_2u_1}+x_{i_1u_1}$. 
The value of $F_{i_1}$ is unchanged; the value of $F_{i_2}$ cannot decrease; and an edge of the support graph disappears.
Otherwise, set
\[
 a_j=\frac{\partial F_{i_j}}{\partial x_{i_ju_j}}>0,\qquad
 b_j=\frac{\partial F_{i_j}}{\partial x_{i_ju_{j-1}}}>0.
\]
Reverse the cycle if necessary to arrange
$\prod_{j=1}^h a_j\ge\prod_{j=1}^h b_j$. Choose $z_1=1$, and for
$2\le j\le h$ set $z_j=b_jz_{j-1}/a_j$. Then
$a_jz_j=b_jz_{j-1}$ for $j\ge2$, and
$a_1z_1\ge b_1z_h$. For each $j$, increase $x_{i_ju_j}$ and decrease
$x_{i_{j+1}u_j}$ by $\varepsilon z_j$. The column sums remain one.
Because the multilinear extension $F_{i}(x)$ is affine in each coordinate $x_{iu}$, the exact change for index $i_j$ is
\[
 \varepsilon(a_jz_j-b_jz_{j-1})
 -\varepsilon^2z_jz_{j-1}
       \frac{\partial^2F_{i_j}}{\partial x_{i_ju_{j-1}}\partial x_{i_ju_{j}}}\ge0.
\]
Increase $\varepsilon$ until a support edge disappears. Repeating
this operation makes the support graph a forest without decreasing any
$F_i$.

Root each nontrivial component of this forest at an index vertex and assign each element vertex $u\in U$ to its parent index vertex. 
If $S_i$ is the final support set of elements $u$ with $x_{iu}>0$, then monotonicity gives $F_i(x_i)\le f_i(S_i)$, while submodularity shows that removing one element costs at most its singleton value, i.e., $f_i(S_i)-f_i(S_i-u)\le f_i(\{u\})$. Thus $f_i(P_i)\ge F_i(x_i)-D_i$.
\end{proof}

For $0\le z\le1$, set
\begin{equation}\label{eq:P}
 \mathcal P(z)=1-(1-z)e^{z-1}.
\end{equation}

\begin{lemma}\label{lem:completion}
Let $H$ be a nonnegatively weighted $r$-uniform multihypergraph of total edge weight $1$.
For distinct vertices $b_1,\ldots,b_r$, put $\sigma=\min_i d_H(b_i)$ and $
 D=\max_{v\notin\{b_1,\ldots,b_r\}}d_H(v)$
with $D=0$ if there are no vertices in $V(H)\setminus \{b_1,\ldots,b_r\}$. Then there is a partition $V(H)=V_1\cup\cdots\cup V_r$ with $b_i\in V_i$ and, for $1\le i\le r$,
\[
 d_H(V_i)\ge \mathcal P(\sigma)-D.
\]
\end{lemma}
\begin{proof}
Write $B=\{b_1,\ldots,b_r\}$, and assign color $i$ to $b_i$.
Given a probability vector $\pi=(\pi_1,\ldots,\pi_r)$, independently assign color $i$ to each vertex outside $B$ with probability $\pi_i$. Then the expected weight missed by color $i$ is $\mathcal D_i(\pi_i)$, where
\[
 \mathcal D_i(p)=\sum_{e:\,b_i\notin e}w(e)(1-p)^{|e\setminus B|}.
\]
Let $g(x)=xe^{-x}$. For every $a\in[0,1]$ and $x\ge0$,
\begin{equation}\label{eq:tangent}
 g(x)\le g(a)+g'(a)(x-a).
\end{equation}
Indeed, $g''(x)=-(2-x)e^{-x}<0$ on $[0,1]$, which implies that \eqref{eq:tangent} holds for $x\in[0,1]$;
for $x\ge1$, since $g'(x)=(1-x)e^{-x}$, we have $g(x)\le g(1)$ and $g'(a)\ge0$, which implies that $g(x)\le g(1)\le g(a)+g'(a)(1-a)\le g(a)+g'(a)(x-a)$.

For a probability vector $\pi=(\pi_1,\ldots,\pi_r)$ and an edge $e$,
let $M_e$ denote the colors absent from $e\cap B$. As $H$ is
$r$-uniform and the prescribed colors are distinct,
$|M_e|=|e\setminus B|$. Put $s_e=\sum_{i\in M_e}\pi_i$.
If $e\nsubseteq B$, then by \eqref{eq:tangent} at $a=s_e$,
\begin{equation}\label{eq:edge-estimate}
 \begin{aligned}
 \sum_{i\in M_e}\pi_i(1-\pi_i)^{|e\setminus B|}
 \le \frac1{|M_e|}\sum_{i\in M_e}g(|e\setminus B|\pi_i)
 \le g(s_e).
 \end{aligned}
\end{equation}
For $e\subseteq B$, both sides are zero. Since $g$ is increasing and concave
on $[0,1]$, Jensen's inequality and \eqref{eq:edge-estimate} give
\begin{equation}\label{eq:weighted-loss}
\begin{aligned}
 \sum_i\pi_i \mathcal D_i(\pi_i)
 &=\sum_{e}\sum_{i\in M_e}w(e)\pi_i(1-\pi_i)^{|e\setminus B|}\ 
 \le\sum_ew(e)g(s_e)\ 
 \le  g\left(\sum_e w(e)s_e\right)\\
 &= g\left(\sum_{i}\sum_{e:b_i\notin e} w(e)\pi_i\right)
 = g\left(1-\sum_i\pi_id_H(b_i)\right)
 \le g(1-\sigma).
\end{aligned}
\end{equation}

If $\sigma=1$, every edge meets each $b_i$, and the assertion is immediate. 
Otherwise, set $T=g(1-\sigma)>0$. For each $i$, let $p_i$ be the least $p\in[0,1]$ such that $\mathcal D_i(p)\le T$; it exists since $\mathcal D_i(1)=0$, and $\mathcal D_i(p_i)=T$ whenever $p_i>0$. 
Let $\pi_i=p_i/(\sum_ip_i)$ if $\sum_ip_i>0$, and $\pi_i=1/r$ if $\sum_ip_i=0$.
If $\sum_i p_i>1$, then $\pi_i<p_i$ and $\mathcal D_i(\pi_i)>T$ whenever $p_i>0$, while $\pi_i=0$ whenever $p_i=0$.
Consequently, $\sum_i\pi_i\mathcal D_i(\pi_i)>T$, contradicting \eqref{eq:weighted-loss}.
Thus we have $\sum_i p_i\le1$, and we obtain a probability vector $\pi$ with $\mathcal D_i(\pi_i)\le T$ for all $i$.

For each $i$, consider the function
\[
 f_i(X)=d_H(\{b_i\}\cup X)-d_H(b_i),
\]
where $X\subseteq V(H)\setminus B$. 
Then each $f_i$ is a normalized, monotone, submodular set function.
Recall that we put all outside vertices into color $i$ with probability $\pi_i$.
Let $R_i$ be the random set consisting of all outside vertices with color $i$.
Then the multilinear extension $F_i$ of $f_i$ satisfies 
$$
F_i(\pi)=\mathbb E[f_i(R_i)]=1-\mathcal D_i(\pi_i)-d_H(b_i)\ge 1-T-d_H(b_i)=\mathcal P(\sigma)-d_H(b_i).
$$
Lemma~\ref{lem:round} gives a partition $P_1\cup\cdots\cup P_r$ such that $f_i(P_i)\ge F_i(\pi_i)-D_i$, where $D_i=\max_{u\in V(H)\setminus B}f_i(\{u\})\le D$.
Let $V_i=P_i\cup\{b_i\}$. Then 
\[
d_H(V_i)=d_H(b_i)+f_i(P_i)\ge d_H(b_i)+ F_i(\pi_i)-D\ge \mathcal P(\sigma)-D.
\]
This completes the proof.
\end{proof}

\begin{corollary}\label{cor:seeds}
For every integer $r\ge2$, let $H$ be a finite nonnegatively weighted $r$-uniform multihypergraph of total edge weight $1$. If there are at most $r$ vertices with degree greater than $t\coloneq\frac{5-\sqrt{13}}6$, then $H$ has a partition $V(H)=V_1\cup\cdots\cup V_r$
such that, for $1\le i\le r$,
\[\ d_H(V_i)\ge c
 \coloneq\frac{11-\sqrt{13}}{18}.
\]
\end{corollary}
\begin{proof}
Suppose at most $r$ vertices have degree greater than $t$.
Let $B$ consist of the $r$ vertices of largest degree.
Let $\delta$ be the maximum degree of a vertex outside $B$. Then $\delta\le t$, and every vertex of $B$ has degree at least $\delta$. 
Use Lemma~\ref{lem:completion} with $D=\delta $ and $\sigma\ge \delta$. There is an $r$-partition such that each class meets at least $\mathcal P(\sigma)-\delta $.
As $\mathcal P(z)$ is increasing and $\mathcal P(z)-z$ is decreasing on $[0,1]$,
we have 
$$\mathcal P(\sigma)-\delta \ge\mathcal P(\delta)-\delta \ge \mathcal P(t)-t\ge c.$$
This completes the proof.
\end{proof}

\section{Hypergraphs with many vertices of large degree}\label{sec:greedy}

In this section, we set
\begin{equation}\label{eq:constants}
c=\frac{11-\sqrt{13}}{18},\qquad
t=3c-1=\frac{5-\sqrt{13}}6,\qquad
\lambda=\frac12-t=\frac{\sqrt{13}-2}{6}.
\end{equation}
These constants satisfy
\begin{equation}\label{eq:relations}
3t^2-5t+1=0,\qquad
t^2+c-2t=0,\qquad
1-2\lambda=2t,
\end{equation}
with $0<t<c<1/2$. 
We work under the assumption that $\Delta(H)<c$.
Put $T=\{v:d_H(v)\ge t\}$. Corollary~\ref{cor:seeds} ensures that a counterexample has at least $r+1$ vertices in $T$.
For $X\subseteq V(H)$ put
\[
 D(X)=\sum_{v\in X}d_H(v)=\sum_{e}w(e)|e\cap X|,\qquad k_X(e)=|e\cap X|.
\]
For $B\subseteq T$, let
\begin{equation}\label{eq:potential}
 \Phi_\lambda(B)=\lambda r D(B)
 +\frac12\sum_e w(e)|e\cap B|^2.
\end{equation}
Choose an $r$-set $B\subseteq T$ maximizing $\Phi_{\lambda}(B)$ over all $r$-subsets of $T$.
For $v\in B$, define
\begin{equation}\label{eq:score}
 a(v)=r\,d_H(v)+\sum_{e\ni v}w(e)|e\cap B|.
\end{equation}
We now describe a greedy algorithm for constructing an $r$-partition. Order \(B=\{v_1,\ldots,v_r\}\) so that \(a(v_1)\le\cdots\le a(v_r)\). Initially, set \(V^{(0)}_i=\{v_i\}\) for each \(i\in[r]\). Select the remaining vertices greedily. Specifically, when selecting \(v_{kr+i}\), choose an unassigned vertex \(u\) that maximizes \(d_H(V^{(k-1)}_i\cup\{u\})\), set \(v_{kr+i}=u\), and update \(V^{(k)}_i\leftarrow V^{(k-1)}_i\cup\{u\}\), where \(k\ge1\) and \(i\in[r]\). Break ties arbitrarily. Let $V_1,\ldots,V_r$ denote the final partition classes.

\begin{lemma}\label{lem:greedy}
Let $U_i=\{v_{r+1},\ldots,v_{r+i-1}\}$ for each $i\le r$ (in particular, $U_1=\varnothing$).
With $V_i$ now denoting the final class, define $Z_i=\sum_{e:\,e\cap V_i=\varnothing}w(e)|e\cap U_i|$.
Then
\begin{equation}\label{eq:greedy-bound}
 (2r-1)d_H(V_i)\ge
 r-D(B)+a(v_i)-d_H(v_i)-Z_i.
\end{equation}
\end{lemma}
\begin{proof}
Fix $i\in[r]$. Note that the edges ultimately missed by $V_i$ have total incidence
$$\sum_{(v,e):\,v\in e,\,e\cap V_i=\varnothing}w(e)=r(1-d_H(V_i)).$$
Now we bound the total incidence of such edges.
Note that the contribution of vertices in $B$ is bounded by
\[
 \sum_{e:\,e\cap V_i=\varnothing}w(e)|e\cap B|\le D(B)-\sum_{e\ni v_i}w(e)|e\cap B|
 =D(B)-a(v_i)+r\,d_H(v_i).
\]
Now consider the vertices $v_j$ with $kr+i+1\le j\le (k+1)r+i-1$ for some $k\ge1$.
The algorithm gives $d_{H}(V^{(k-1)}_i\cup\{v_{kr+i}\})\ge d_{H}(V^{(k-1)}_i\cup\{v_{j}\})$ for each $j$.
Thus, the contribution of each such $v_j$ is at most
$$\sum_{(v_j,e):\,v_j\in e,\,e\cap V_i=\varnothing}w(e)\le d_{H}(V^{(k-1)}_i\cup\{v_{j}\})-d_{H}(V^{(k-1)}_i)\le d_{H}(V^{(k)}_i)-d_{H}(V^{(k-1)}_i).$$
Therefore, the contribution of the vertices $v_j$ with $kr+i+1\le j\le (k+1)r+i-1$ is at most $(r-1)\left(d_{H}(V^{(k)}_i)-d_{H}(V^{(k-1)}_i)\right)$. Summing over $k\ge1$, we find that the total contribution of these later vertices is at most
\begin{align*}
    (r-1)\sum_{k\ge1}\left(d_{H}(V^{(k)}_i)-d_{H}(V^{(k-1)}_i)\right)= (r-1)\left(d_H(V_i)-d_H(v_i)\right).
\end{align*}
Moreover, the vertices in $U_i$ contribute $Z_i$. Consequently,
\begin{align*}
 r(1-d_H(V_i))
 &\le D(B)-a(v_i)+r\,d_H(v_i)+Z_i+(r-1)(d_H(V_i)-d_H(v_i))\\
 &=D(B)-a(v_i)+d_H(v_i)+Z_i
 +(r-1)d_H(V_i),
\end{align*}
which completes the proof.
\end{proof}

Recall that $V_1,\ldots,V_r$ are the partition classes obtained by the greedy algorithm. 
Fix an index $i$ such that $d_H(V_i)<c$.
If no such index exists, the partition already satisfies Theorem~\ref{thm:main}.
Write
\[
 P=\{v_1,\ldots,v_i\},\quad R=B\setminus P,\quad
 \alpha=\frac ir,\quad
 A=\frac{a(v_i)}r,\quad
 \overline A=\frac1{ir}\sum_{x\in P}a(x).
\]
Let $U\subseteq T\setminus B$ consist of $\ell\le i$
vertices, and set $\theta=\ell/i$. 
Sample an edge $e$ of $H$ with $\Pr[e=f]=w(f)$ for every $f\in E(H)$.
Write
\begin{equation}\label{eq:variables}
 X=\frac{|e\cap P|}i,\qquad
 Y=\begin{cases}|e\cap U|/\ell,&\ell>0,\\0,&\ell=0,\end{cases}
 \qquad W=\frac{|e\cap R|}r.
\end{equation}
Thus $0\le X,Y\le1$ and $0\le W\le1-\alpha$.

\begin{claim}\label{clm:coupled}
With the above notation,
\begin{equation}\label{eq:identity}
 A\ge\overline A=\mathbb E X+\alpha\mathbb E  X^2+\mathbb E (XW),
\end{equation}
\begin{equation}\label{eq:exchange}
\begin{aligned}
 \overline A\ge(1-\lambda)\mathbb E X+\lambda\mathbb E Y
 +\frac{\alpha\theta}{2}\mathbb E(X^2+Y^2)+\alpha(1-\theta)\mathbb E(XY)+\mathbb E(WY),
\end{aligned}
\end{equation}
and, under $\Delta(H)<c$,
\begin{equation}\label{eq:suffix}
 \mathbb E  W\le(1-\alpha)c,\qquad
 \mathbb E (XW)\le\frac{1-\alpha}{2}(\mathbb E  X^2+c).
\end{equation}
\end{claim}
\begin{proof}
The ordering of $B$ gives $A\ge\overline A$. Summing \eqref{eq:score} over $P$ proves \eqref{eq:identity}.

If $\ell=0$, \eqref{eq:exchange} follows directly from $\overline A\ge \mathbb EX$. 
If $\ell>0$, select uniformly an $\ell$-subset $P'\subseteq P$.
The maximality of $B$ in \eqref{eq:potential} implies that $\Phi_{\lambda}(B)-\Phi_{\lambda}(B\setminus P')\ge\Phi_{\lambda}(B\setminus P'\cup U)-\Phi_{\lambda}(B\setminus P')$.
Note that the loss from deleting $P'$ equals
\[
 \Phi_{\lambda}(B)-\Phi_{\lambda}(B\setminus P')=\sum_{x\in P'}a(x)-(1-\lambda)r D(P')
 -\frac12\sum_e w(e)|e\cap P'|^2,
\]
and the gain from adding $U$ to $B\setminus P'$ equals
\[
 \Phi_{\lambda}(B\setminus P'\cup U)-\Phi_{\lambda}(B\setminus P')
 =\lambda rD(U)+\frac12\sum_e w(e)|e\cap U|^2
 +\sum_e w(e)|e\cap (B\setminus P')||e\cap U|.
\]
Moreover, note that
 $\mathbb E_{P'}|e\cap P'|^2\ge\theta^2|e\cap P|^2$ and 
 $\mathbb E_{P'}|e\cap(B\setminus P')|=|e\cap R|+(1-\theta)|e\cap P|$.
Hence we have
$$\ell r\overline A\ge(1-\lambda)\ell r \mathbb EX+\lambda\ell r\mathbb EY
 +\frac{\ell^2}2\mathbb EX^2+\frac{\ell^2}{2}\mathbb E(Y^2)
 +\ell\, \mathbb E((rW+(1-\theta)iX)Y).$$
Dividing by $\ell r$ yields \eqref{eq:exchange}. 

Finally, $|R|=(1-\alpha)r$ and $\Delta(H)\le c$ imply
$\mathbb E W=D(R)/r\le(1-\alpha)c$. For $\alpha<1$, the inequality
$2XW\le(1-\alpha)X^2+W^2/(1-\alpha)$ gives
\[
 2\mathbb E (XW)\le(1-\alpha)\mathbb E  X^2+
 \frac{\mathbb E  W^2}{1-\alpha}
 \le(1-\alpha)(\mathbb E  X^2+c),
\]
since $W^2\le(1-\alpha)W$. For $\alpha=1$, we have $W=0$.
\end{proof}

Recall that we randomly choose an edge $e$ with probability $w(e)$. Let $M$ be the event that $e\cap V_i=\varnothing$. Define $y=\mathbb E (Y\mathbf 1_M)$ and
\begin{equation*}
 Z=t+\alpha(\mathbb E X+t-c)+\alpha\theta(y-t).
\end{equation*}
In the following two subsections, we derive upper and lower bounds for $A-Z$ to obtain a contradiction.
\subsection{Upper bound}
We first establish the upper bound supplied by the greedy algorithm.

\begin{lemma}\label{lem:upper}
Assume $\Delta(H)\le c$ and carry out the construction of this section. 
Fix $i\in\{1,\ldots,r\}$, and put $q=d_H(V_i),\ U=U_i\cap T,\ \ell=|U|$.
Use the notation of \eqref{eq:variables} with this choice of $U$, so that
$\alpha=i/r$ and $\theta=\ell/i$. Let $M=\{e:e\cap V_i=\varnothing\}$ and $y=\mathbb E (Y\mathbf 1_M),$ where the edge is sampled according to its weight. 
For $Z=t+\alpha(\mathbb E X+t-c)+\alpha\theta(y-t)$, we have
\begin{equation}\label{eq:upper-strict}
 A-Z<-2(c-q).
\end{equation}
\end{lemma}

\begin{proof}
Since $H$ has total edge weight $1$, we have $\sum_{v\in V(H)}d_H(v)=r$.
Since every vertex has degree at most $c<1/2$, more than $2r$ vertices have positive degree, which means that at least $r$ vertices lie outside $B$.
Therefore, we have $|U_i|=i-1$.
Recall that $Z_i=\sum_{e\in M}w(e)|e\cap U_i|$.
Among those $i-1$ vertices of $U_i$, there are exactly $\ell$ vertices lying in $U$; they contribute
\[
 \sum_{e\in M}w(e)|e\cap U|=\ell y
\]
to $Z_i$ (with both sides zero when $\ell=0$). All other early vertices
lie outside $T$, so each has degree strictly below $t$. It
follows that
\begin{equation}\label{eq:upper-Zi}
 \frac{Z_i}{r}\le
 \alpha\theta y+\left(\alpha-\alpha\theta-\frac1r\right)t.
\end{equation}
Furthermore, since $D(P)=i\mathbb EX$ and $d_H(v)\le c$ for $v\in R$,
\begin{equation}\label{eq:upper-SB}
 \frac{D(B)}r=\frac{D(P)+D(R)}r\le\alpha \mathbb EX+(1-\alpha)c.
\end{equation}
By Lemma~\ref{lem:greedy} and \eqref{eq:greedy-bound},
\[
 A\le \left(2-\frac 1r\right)q-1+\frac{D(B)}r+\frac {d_H(v_i)}r+\frac{Z_i}r.
\]
Substituting \eqref{eq:upper-Zi} and \eqref{eq:upper-SB} yields
\[
\begin{aligned}
 A&\le 2q-1+\alpha \mathbb EX+(1-\alpha)c+\alpha\theta y
      +(\alpha-\alpha\theta)t-\frac{t+q-d_H(v_i)}{r}\\
  &=t+\alpha(\mathbb EX+t-c)+\alpha\theta(y-t)-2(c-q)
      -\frac{t+q-d_H(v_i)}{r}\\
  &=Z-2(c-q)-\frac{t+q-d_H(v_i)}{r},
\end{aligned}
\]
where the second equality uses $t=3c-1$. Since $d_H(v_i)\le q$ and $t>0$, we have $t+q-d_H(v_i)>0$, proving \eqref{eq:upper-strict}.
%
\end{proof}

\subsection{Lower bound}

We next establish a lower bound for $A-Z$.

\begin{lemma}\label{lem:moments}
Let $c,t,\lambda$ be as in \eqref{eq:constants}. Suppose $0<\alpha\le1$, $0\le\theta\le1$, $t\le x\le c$, $0\le y\le c$, $0\le q<c$,
and let $X,Y,W$ be nonnegative random variables satisfying $\mathbb E X=x$, \eqref{eq:identity}, \eqref{eq:exchange} and \eqref{eq:suffix}.
Let $M$ be an event of probability $1-q$ such that $y=\mathbb E(Y\mathbf 1_M)$.
With $Z=t+\alpha(x+t-c)+\alpha\theta(y-t)$, we have
\begin{equation*}
 A\ge Z-\frac{c-q}{3}.
\end{equation*}
\end{lemma}

\begin{proof}
Suppose first that $y\le t$. Jensen's inequality applied to
\eqref{eq:identity} gives $A\ge x+\alpha x^2$. As
$\theta(y-t)\le0$, $t^2+c-2t=0$ and $x=D(P)/i\ge t$, we obtain
\[
\begin{aligned}
 A-Z&\ge x-t+\alpha(x^2-x+c-t)\\
 &=(1-\alpha)(x-t)+\alpha(x^2-t^2)\ge0.
\end{aligned}
\]
We may therefore assume $y>t$, so both $x$ and $y$ lie in $[t,c]$.
For $\mu\ge0$, \eqref{eq:identity} and \eqref{eq:exchange} give
\begin{align*}
    A\ge(1-\mu)\overline A+\mu\overline A\ge{}&(1-\mu\lambda)\mathbb EX+\mu\lambda \mathbb EY
 +\alpha\Bigl(1-\mu+\frac{\mu\theta}{2}\Bigr)\mathbb E  X^2\\
 &+\frac{\mu\alpha\theta}{2}\mathbb E  Y^2
 +\mu\alpha(1-\theta)\mathbb E (XY)
 +(1-\mu)\mathbb E (XW)+\mu \mathbb E(WY).
\end{align*}

We first consider $0\le\theta\le9/20$.
Let $\mu=2\theta$. Then we have
\begin{align*}
    A\ge(1-2\lambda)\mathbb EX+2\lambda\mathbb E[(1-\theta)X+\theta Y]+\alpha\mathbb E\bigl[((1-\theta)X+\theta Y)^2\bigr]+\mathbb E\{W[(1-2\theta)X+2\theta Y]\}.
\end{align*}
Note that $(1-2\theta)X+2\theta Y\ge \min\{X,Y\}\ge0$. Put
$z=(1-\theta)x+\theta y$. Jensen's inequality, together with $\mathbb EX=x$, $\mathbb EY\ge y$, gives
\[
 A\ge (1-2\lambda)x+2\lambda z+\alpha z^2.
\]
Moreover, with $z=(1-\theta)x+\theta y$,
\begin{align*}
    Z=t+\alpha(x+t-c)+\alpha\theta(y-t)=t+\alpha(t-c+z+\theta (x-t)).
\end{align*}
So
\begin{align*}
    A-Z\ge (1-2\lambda)x+2\lambda z-t+\alpha(z^2-t+c-z-\theta (x-t)).
\end{align*}
Apply \eqref{eq:relations}. Using $2t+2\lambda=1$, we have
\begin{align*}
    (1-2\lambda)x+2\lambda z-t=2tx+(1-2t)z-t=2t(x-t)+(1-2t)(z-t).
\end{align*}
Using $c=2t-t^2$, we have
\begin{align*}
    z^2-t+c-z=z^2-z+t-t^2=(z-t)^2-(1-2t)(z-t).
\end{align*}
Combining the two identities above, together with $2t+2\lambda=1$, gives
\begin{equation*}
\begin{aligned}
 A-Z
 &\ge2t(x-t)+(1-2t)(z-t)+\alpha[(z-t)^2-(1-2t)(z-t)-\theta (x-t)]\\
 &={}(1-\alpha)
 \bigl[2t(x-t)+(1-2t)(z-t)\bigr]+\alpha\bigl[(z-t)^2+(2t-\theta)(x-t)\bigr].
\end{aligned}
\end{equation*}
Since $x,z\ge t$ and $1\ge 2t\ge9/20$, we have
\begin{align}\label{eq:early-square}
    A-Z\ge0\ge -\frac{c-q}{3}.
\end{align}

Now we consider $9/20\le\theta\le1$.
Use $\mathbb E Y\ge y$ and $\mathbb E(WY)\ge0$. We have
\begin{equation*}
\begin{aligned}
 A\ge{}&(1-\mu\lambda)x+\mu\lambda y
 +\alpha\Bigl(1-\mu+\frac{\mu\theta}{2}\Bigr)\mathbb E  X^2\\
 &+\frac{\mu\alpha\theta}{2}\mathbb E  Y^2
 +\mu\alpha(1-\theta)\mathbb E (XY)
 +(1-\mu)\mathbb E (XW).
\end{aligned}
\end{equation*}
Let $\beta=(\mu-1)_+$ and define
\begin{equation*}
\begin{aligned}
 a=\alpha\left(1-\mu+\frac{\mu\theta}{2}\right)
       -\frac{\beta(1-\alpha)}2,
       \qquad
 b=\frac{\mu\alpha(1-\theta)}2,
 \qquad 
 d=\frac{\mu\alpha\theta}2.
\end{aligned}
\end{equation*}
The bound on $\mathbb E (XW)$ in \eqref{eq:suffix} now implies
\begin{equation}\label{eq:quadratic}
 A\ge(1-\mu\lambda)x+\mu\lambda y
      -\frac{\beta(1-\alpha)c}{2}
      +\mathbb E (aX^2+2bXY+dY^2).
\end{equation}

Choose $\mu$ according to the following five rectangles:
\begin{center}
\begin{tabular}{@{}ccc@{}}
\toprule
$\theta$-interval&$\alpha$-interval&$\mu$\\
\midrule
$[9/20,3/5]$&$(0,1]$&$9/10$\\
$(3/5,4/5]$&$(0,1]$&$1$\\
$(4/5,1]$&$(0,7/10]$&$1$\\
$(4/5,19/20]$&$(7/10,1]$&$11/10$\\
$(19/20,1]$&$(7/10,1]$&$6/5$\\
\bottomrule
\end{tabular}
\end{center}
Appendix~\ref{app:certificates_1} shows that, in every rectangle, $a>0$ and $ad-b^2\ge0$.
Note that 
\begin{align*}
    \mathbb E(aX^2+2bXY+dY^2)=a\mathbb E\left[\left(X+\frac{b}{a}Y\right)^2\right]+\left(d-\frac{b^2}{a}\right)\mathbb EY^2.
\end{align*}
By the Cauchy--Schwarz inequality, $\mathbb EY^2\mathbb E(\mathbf 1_M)\ge [\mathbb E(Y\mathbf1_M)]^2=y^2$. 
Thus, by $c\ge q$,
\begin{align*}
    \left(d-\frac{b^2}{a}\right)\mathbb EY^2
    \ge\left(d-\frac{b^2}{a}\right)\frac{y^2}{1-q}
    \ge\left(d-\frac{b^2}{a}\right)\frac{y^2}{1-c}-\frac{dy^2(c-q)}{(1-c)^2}.
\end{align*}
Since $y\le \mathbb EY=D(U)/\ell\le c$ and $d=\mu\alpha\theta/2\le\mu/2$,
\begin{align*}
    \frac{dy^2(c-q)}{(1-c)^2}\le \frac{\mu c^2}{2(1-c)^2}(c-q)\le\frac{c-q}{3},
\end{align*}
where the last inequality holds because our choice of $\mu$ is at most $6/5$ in each case and $c^2/(1-c)^2\le5/9$.
Moreover, by Jensen's inequality, together with $\mathbb EX=x\ge0$ and $\mathbb EY\ge y\ge0$, we have $\mathbb E\left[\left(X+\frac{b}{a}Y\right)^2\right]\ge\left(x+\frac{b}{a}y\right)^2$.
Combining these inequalities, together with $ad-b^2\ge0$, gives
\begin{equation}\label{eq:cs}
\begin{aligned}
 \mathbb E (aX^2+2bXY+dY^2)
 &\ge a\left(x+\frac ba y\right)^2
 +\left(d-\frac{b^2}{a}\right)\frac{y^2}{1-c}-\frac{c-q}{3}.
\end{aligned}
\end{equation}
Put
\begin{equation*}
\begin{aligned}
 C=\frac{d-cb^2/a}{1-c},\
 L=1-\mu\lambda-\alpha,\
 N=\mu\lambda-\alpha\theta,\
 K=-t-\alpha(t-c)+\alpha\theta t
          -\frac{\beta(1-\alpha)c}{2}.
\end{aligned}
\end{equation*}
Consequently, \eqref{eq:quadratic} and \eqref{eq:cs} imply $A-Z\ge F_{\alpha,\theta,\mu}(x,y)-\frac{c-q}{3}$, where
\begin{equation*}
 F_{\alpha,\theta,\mu}(x,y)
 =ax^2+2bxy+Cy^2+Lx+Ny+K.
\end{equation*}
Appendix~\ref{app:certificates_2} shows that $F_{\alpha,\theta,\mu}(x,y)\ge0$ for $x,y\in[t,c]$ on the five rectangles. So
\begin{align*}
    A-Z\ge F_{\alpha,\theta,\mu}(x,y)-\frac{c-q}{3}\ge-\frac{c-q}{3}.
\end{align*}
This together with \eqref{eq:early-square} completes the proof.
\end{proof}

\section{Proof of the partition theorem}\label{sec:proof}

\begin{proof}[Proof of Theorem~\ref{thm:main}]
We proceed by induction on $r$. For $r=2$, this follows from the local optimality of a maximum weighted cut. 

Suppose $r\ge3$ and the result is known for $r-1$. If there is a vertex $v$ with $d_H(v)\ge c$, take $\{v\}$ as one class. 
For each edge $e\in E(H)$, choose an $(r-1)$-subset $f_e\subseteq e\setminus\{v\}$, and define a weighted $(r-1)$-uniform hypergraph $H'$ by $w'(f)=\sum_{e\in E(H):\,f_e=f}w(e)$.
Then $H'$ has total edge weight $1$. 
By induction, its vertices admit an $(r-1)$-partition in which every class meets edges of $H'$ of total weight at least $c$, and hence every class meets edges of $H$ of total weight at least $c$. 
We may therefore assume that $\Delta(H)<c$.

Suppose, for a contradiction, that there is no required partition. 
Corollary~\ref{cor:seeds} implies that there are at least $r+1$ vertices of degree at least $t$.
Run the cyclic greedy allocation and suppose that some class $V_i$ meets edges of total weight less than $c$, i.e., $q=d_H(V_i)<c$.
Specialize the variables as in Lemma~\ref{lem:upper}, taking $U=U_i\cap T$ and $y=\mathbb E (Y\mathbf 1_M)$ for the event $M$ of edges missed by $V_i$. 
Then Lemma~\ref{lem:upper} gives
\[
 A-Z<-2(c-q).
\]
Moreover, Lemma~\ref{lem:moments} gives
\[
 A-Z\ge-\frac{c-q}{3}.
\]
Combining these two bounds, we obtain $-2(c-q)>-\frac{c-q}{3}$, which means $q>c$, contradicting the assumption $q<c$. This completes the proof.
\end{proof}

\section*{Declaration on the Use of AI}

The authors used ChatGPT 5.6 to assist in discussing proof strategies, checking proofs, and improving exposition.

\newpage

\appendix
\section{Exact calculations in Lemma~\ref{lem:moments}}\label{app:certificates}

This appendix records finite, exact arithmetic that verifies the
corner calculations in Lemma~\ref{lem:moments}. The calculations
use exact arithmetic and the identities
\begin{equation}\label{eq:t-rational}
 t=\frac{5-\sqrt{13}}{6},
 \qquad 3t^2-5t+1=0.
\end{equation}

Recall that we partition the rectangle $\theta\in[9/20,1],\alpha\in(0,1]$ into five smaller rectangles and choose the value of $\mu$ for each small rectangle in the proof of Lemma~\ref{lem:moments}:
\begin{center}
\begin{tabular}{@{}ccc@{}}
\toprule
$\theta$-interval&$\alpha$-interval&$\mu$\\
\midrule
$[9/20,3/5]$&$(0,1]$&$9/10$\\
$[3/5,4/5]$&$(0,1]$&$1$\\
$[4/5,1]$&$(0,7/10]$&$1$\\
$[4/5,19/20]$&$(7/10,1]$&$11/10$\\
$[19/20,1]$&$(7/10,1]$&$6/5$\\
\bottomrule
\end{tabular}
\end{center}

\subsection{Positive semidefiniteness of the quadratic form}\label{app:certificates_1}
For the first three small rectangles, that is, $9/20\le\theta\le4/5,\,0<\alpha\le1$ and $4/5\le\theta\le1,\,0<\alpha\le7/10$, we have $\mu\le1$, $a>0$ and 
\begin{equation}\label{eq:det}
 ad-b^2=\frac{\alpha^2\mu(2\theta-\mu)}4\ge0.
\end{equation}
For the remaining case $4/5\le \theta\le1,\,7/10\le\alpha\le1$, we have $a=\alpha(\frac{11}{20}\theta-\frac{1}{20})-\frac1{20}>0$ when $\mu=11/10$ and $a=\alpha(\frac35\theta-\frac{1}{10})-\frac{1}{10}>0$ when $\mu=6/5$.
Note that, for any small rectangle $\alpha_0\le\alpha\le\alpha_1$, $\theta_0\le\theta\le\theta_1$, each matrix $Q(\alpha,\theta)=\bigl(\begin{smallmatrix}a&b\\b&d\end{smallmatrix}\bigr)$ can be represented by
\begin{align*}
    Q(\alpha,\theta)=(1-s_1)(1-s_2)Q(\alpha_0,\theta_0)+(1-s_1)s_2 Q(\alpha_0,\theta_1)+s_1(1-s_2)Q(\alpha_1,\theta_0)+s_1s_2 Q(\alpha_1,\theta_1),
\end{align*}
where $s_1=\frac{\alpha-\alpha_0}{\alpha_1-\alpha_0}$ and $s_2=\frac{\theta-\theta_0}{\theta_1-\theta_0}$. Thus $Q(\alpha,\theta)$ is positive semidefinite if the four matrices $Q(\alpha_i,\theta_j)$, $i,j\in\{0,1\}$, are positive semidefinite.
The following exact corner values give strict positivity throughout both rectangles:
\begin{center}
\begin{tabular}{@{}ccccc@{}}
\toprule
$\mu$&$\alpha$&$\theta$&$a$&$ad-b^2$\\
\midrule
$11/10$&$7/10$&$4/5$&$223/1000$&$12551/200000$\\
$11/10$&$7/10$&$19/20$&$1123/4000$&$81851/800000$\\
$11/10$&$1$&$4/5$&$17/50$&$11/80$\\
$11/10$&$1$&$19/20$&$169/400$&$11/50$\\
$6/5$&$7/10$&$19/20$&$229/1000$&$9093/100000$\\
$6/5$&$7/10$&$1$&$1/4$&$21/200$\\
$6/5$&$1$&$19/20$&$37/100$&$21/100$\\
$6/5$&$1$&$1$&$2/5$&$6/25$\\
\bottomrule
\end{tabular}
\end{center}

\subsection{A lower bound for $F_{\alpha,\theta,\mu}$}\label{app:certificates_2}
Recall that we defined
\[
\begin{aligned}
F_{\alpha,\theta,\mu}(x,y)
={}&ax^2+2bxy+\frac{d}{1-c}y^2
+Lx+Ny+K
-\frac{cb^2}{(1-c)a}y^2,
\end{aligned}
\]
in the proof of Lemma~\ref{lem:moments}. 
Set $F_{0,\theta,\mu}(x,y)=(1-\mu\lambda)x+\mu\lambda y-t$.
The coefficients $a,b,d,L,N,K$ are affine in each of $\alpha$ and $\theta$ separately.
Thus, for $z\in\{\alpha,\theta\}$,
\[
\frac{\partial^2 F_{\alpha,\theta,\mu}}
{\partial z^2}
=
-\frac{2cy^2}{(1-c)a^3}
(ab_z-ba_z)^2\le0.
\]
Thus $F_{\alpha,\theta,\mu}(x,y)$ is concave in each of $\alpha$ and $\theta$ separately.
Let $G(\alpha,\theta)=F_{\alpha,\theta,\mu}(x,y)$.
Then, in each small rectangle $\alpha_0 \le \alpha \le \alpha_1, \theta_0\le \theta \le \theta_1$, for fixed $x,y,\mu$, we have 
\begin{align*}
    G(\alpha,\theta)\ge(1-s_1)(1-s_2)G(\alpha_0,\theta_0)+(1-s_1)s_2 G(\alpha_0,\theta_1)+s_1(1-s_2)G(\alpha_1,\theta_0)+s_1s_2 G(\alpha_1,\theta_1),
\end{align*}
where $s_1=\frac{\alpha-\alpha_0}{\alpha_1-\alpha_0}$ and $s_2=\frac{\theta-\theta_0}{\theta_1-\theta_0}$. 
This implies that $F_{\alpha,\theta,\mu}(x,y)\ge0$ if $F_{\alpha_i,\theta_j,\mu}(x,y)\ge0$ for $i,j\in\{0,1\}$.

We now prove that $F_{\alpha,\theta,\mu}(x,y)\ge0$ for $x,y\in[t,c]$ at every corner of the five rectangles. 
For $\alpha=0$, note that $F_{0,\theta,\mu}(x,y)=(1-\mu\lambda)x+\mu\lambda y-t\ge(1-\mu\lambda)t+\mu\lambda t-t=0$.
For $(\alpha,\theta,\mu)=(1,\frac{9}{20},\frac{9}{10})$, recall that, in the proof of Lemma~\ref{lem:moments}, we have proved that 
\begin{align*}
    F_{1,9/20,9/10}(x,y)=\frac{[11(x-t)+9(y-t)]^2}{400}+\left(2t-\frac{9}{20}\right)(x-t)\ge0.
\end{align*}
For the remaining cases, fix $\alpha,\theta,\mu$. We use two lower bounds for $F_{\alpha,\theta,\mu}(x,y)$.

For a row of \emph{type I}, we have $F_x=2ax+2by+L\ge2(a+b)t+L$. 
If $2(a+b)t+L>0$, then $F$ is increasing in $x$.
Thus, $F(x,y)\ge F(t,y)=Cy^2+(2bt+N)y+K+at^2+Lt$.
Thus we have the first lower bound
\begin{align*}
    \Lambda^{\rm I}(\alpha,\theta,\mu)=K+at^2+Lt-\frac{(N+2bt)^2}{4C}.
\end{align*}
For a row of \emph{type II}, note that $\bigl(\begin{smallmatrix}a&b\\b&C\end{smallmatrix}\bigr)$ is positive definite since $aC-b^2>0$. Thus, over $(x,y)\in\mathbb R^2$, $F(x,y)$ has the lower bound
\begin{equation}\label{eq:type-II}
 \Lambda^{\rm II}=
 K-\frac{CL^2-2bLN+aN^2}{4(aC-b^2)}.
\end{equation}
These formulas are valid because $C>0$ and $aC-b^2>0$ at the
listed endpoints. For the type-I rows numbered
$1,4,5,6,9,10$ below, the values of $2(a+b)t+L$ are
\[
 -\frac9{20}+2t,\quad -\frac15+\frac{17}{10}t,\quad
 -\frac14+\frac{17}{10}t,\quad
 -\frac14+\frac{17}{10}t,\quad
 -\frac3{10}+\frac{17}{10}t,\quad
 -\frac3{10}+\frac{17}{10}t;
\]
all are positive by \eqref{eq:t-rational}.

\begin{center}
\renewcommand{\arraystretch}{1.13}
\begin{tabular}{@{}rccccr@{}}
\toprule
No.&$\mu$&$\alpha$&$\theta$&Type&$10^6\Lambda>$\\
\midrule
1&$9/10$&$1$&$3/5$&I&$1672$\\
2&$1$&$1$&$3/5$&II&$4793$\\
3&$1$&$1$&$4/5$&II&$2038$\\
4&$1$&$7/10$&$1$&I&$2900$\\
5&$11/10$&$7/10$&$4/5$&I&$3615$\\
6&$11/10$&$7/10$&$19/20$&I&$4567$\\
7&$11/10$&$1$&$4/5$&II&$6851$\\
8&$11/10$&$1$&$19/20$&II&$904$\\
9&$6/5$&$7/10$&$19/20$&I&$861$\\
10&$6/5$&$7/10$&$1$&I&$1079$\\
11&$6/5$&$1$&$19/20$&II&$2781$\\
12&$6/5$&$1$&$1$&II&$800$\\
\bottomrule
\end{tabular}
\end{center}
For the remaining corner $(\alpha,\theta,\mu)=(7/10,4/5,1)$, the inequality $F_{7/10,4/5,1}(x,y)\ge0$ follows from concavity in $\alpha$, since $F_{1,4/5,1}(x,y)\ge0$ and $F_{0,4/5,1}(x,y)\ge0$.
Therefore, $F_{\alpha,\theta,\mu}(x,y)\ge0$ for $x,y\in[t,c]$ throughout all five rectangles.

\end{document}